\documentclass[12pt, english]{article}
\usepackage[T2A]{fontenc}
\usepackage[cp1251]{inputenc}
\usepackage{babel}

\usepackage[centertags]{amsmath}
\usepackage{amsfonts}
\usepackage{amssymb}
\usepackage{amsthm}
\usepackage{newlfont}
\usepackage[active]{srcltx}

\newtheorem{theorem}{\bf Theorem}
\newtheorem{lemma}{\bf Lemma}
\newtheorem{remark}{\bf Remark}
\newtheorem{corollary}{\bf Corollary}

\begin{document}

\title{On lower and upper bounds for probabilities of unions 
and the Borel--Cantelli lemma}

\author{Andrei N. Frolov
\\ Dept. of Mathematics and Mechanics
\\ St.~Petersburg State University
\\ St. Petersburg, Russia
\\ E-mail address: Andrei.Frolov@pobox.spbu.ru}

\maketitle

{\abstract{
We obtain new lower and upper bounds for probabilities of unions of events.
These bounds are sharp. They are stronger than earlier ones.
General bounds maybe applied in arbitrary measurable spaces.
We have improved the method that has been introduced in previous papers.
We derive new generalizations of the first and second part of the
Borel--Cantelli lemma.
 }


}

\section{Introduction}

Let $(\Omega,\mathcal{F},\mathbf{P})$ be a probability space
and $\{A_n\}$ be a sequence of events. Put $U_n=\bigcup\limits_{i=1}^n A_i$
and $\xi_n= I_{A_1}+I_{A_2}+\cdots+I_{A_n}$ for $n=1,2,\ldots$,
where $I_B$ denotes the indicator of event $B$.

Bounds for $\mathbf{P} (U_n)$ play an important role in
probability and statistics. Many of them are based on
moments $\alpha_k(n) = \mathbf{E} \xi_n^k$, $k=1, 2, \ldots$,
of random variable $\xi_n$.
For example, Bonferroni type inequalities use binomial moments
of $\xi_n$, but one can write them in terms of $\alpha_k(n)$ as well.

Chung and Erd\H{o}s (1952) have derived the most simple and
applicable lower bound for $\mathbf{P} (U_n)$ of this type that is
\begin{eqnarray*}
\mathbf{P}(U_n) \geqslant \frac{\alpha_1^2(n)}{\alpha_2(n)}.
\end{eqnarray*}
It follows from the Cauchy--Buniakowski inequality.
It is very convenient because of
$$ \alpha_1(n) = \sum\limits_{i=1}^n \mathbf{P} (A_i), \quad
\alpha_2(n) = \sum\limits_{i=1}^n \sum\limits_{j=1}^n \mathbf{P} (A_i A_j).
$$

Various generalizations of the Chung--Erd\H{o}s inequality
were obtained in Dawson and Sankoff (1967), Gallot (1966), Kounias (1968),
Kwerel (1975), Boros and Pr\'{e}kopa (1989), Galambos and Simonelli (1996),
de Caen (1997), Kuai, Alajaji and Takahara (2000),
Pr\'{e}kopa (2009), Frolov (2012) and references therein.

For example, in Kwerel (1975) and Boros and Pr\'{e}kopa (1989), one can find
upper and lower bounds for $\mathbf{P}(U_n)$ that are based on
$\alpha_k(n)$ for $1\leqslant k \leqslant 3$ and $\alpha_k(n)$
for $1\leqslant k \leqslant 4$. The lower bounds are stronger
than the Chung--Erd\H{o}s inequality. They are simple enough in applications
as well. Indeed, for every $k \leqslant n$, moments $\alpha_k(n)$ are sums
of probabilities of intersections of $k$ events from
$A_1, A_2, \ldots, A_n$. Note that a precision of bounds increases when
numbers of used moments enlarge.

From the other hand, making use of the H\"{o}lder inequality,
one can conclude that
\begin{eqnarray*}
\mathbf{P}(U_n)  \geqslant \lim\limits_{p\rightarrow 1+}
\left(\frac{(\mathbf{E} \xi_n)^p}{\mathbf{E} \xi_n^p }\right)^{1/(p-1)}
 \geqslant
\left(\frac{(\mathbf{E} \xi_n)^p}{\mathbf{E} \xi_n^p }\right)^{q/p},
\end{eqnarray*}
where $p>1$ and $1/p+1/q =1$.
It follows that applications of $L_p$-norms of $\xi_n/\mathbf{E} \xi_n$ with $p>2$
do not give lower bounds stronger than the Chung--Erd\H{o}s inequality.
In particular, one can not derive better bounds, using
$\alpha_2(n)$ and $\alpha_3(n)$ instead of $\alpha_1(n)$ and $\alpha_2(n)$,
for instance. Stronger bounds have to involve moments of smaller orders and
then we have to use
moments of non-integer orders. Of course, these bounds are more complicated
in calculations. One can find such lower bounds in Frolov (2012)
and upper bounds in Frolov (2014).
Note that in Frolov (2012), one part of bounds is obtained
by applications of the Cauchy--Buniakowski and  H\"{o}lder inequalities
and another one is proved by a method which we improve in this paper.
One can find a discussion on relationship of these bounds in there.

In this paper, we improve the method from Frolov (2012) and
derive new upper and lower bound for $\mathbf{P} (U_n)$.
Some of them include quantities similar to $\alpha_k(n)$,
and may be calculated relatively simple.

Every new bound for $\mathbf{P} (U_n)$ may be used to
obtain generalizations of the Borel-Cantelli lemma,
the classical variant of which is as follows.

{\bf The Borel-Cantelli lemma.}
{\it 1) If series $\sum\limits_{n=1}^\infty \mathbf{P}(A_n)$ converges,
then $\mathbf{P}(A_n \: i.o.)=0$, where
$$ \{A_n \: i.o.\} = \bigcap\limits_{n=1}^\infty\bigcup\limits_{k=n}^\infty A_k.
$$

2) If $\{A_n\}$ are independent and series
$\sum\limits_{n=1}^\infty \mathbf{P}(A_n)$ diverges,
then $\mathbf{P}(A_n \: i.o.)=1$.
}


The first part of the Borel-Cantelli lemma works in many situations, but
the assumption of independence in the second one is restrictive.
Therefore, main attention was paid to generalization of the second part.

The most important generalization of the second part of the Borel-Cantelli lemma
has been obtained by Erd\H{o}s and R\'{e}nyi (1959).
They proved that $\mathbf{P}(A_n \: i.o.)=1$ provided $L=1$,
where
$$ L = \liminf\limits_{n \rightarrow \infty} \frac{\alpha_2(n)}{\alpha_1^2(n)}.
$$
This result implies in particular that the independence maybe relaxed to
pairwise independence in the second part.
Of course, the proof is based on the Chung--Erd\H{o}s inequality.
Kochen and Stone (1964) and Spitzer (1964) have proved that
$\mathbf{P}(A_n \: i.o.)\geqslant 1/L$.
Further generalizations of the second part of the Borel--Cantelli lemma
have been obtained in Kounias (1968), M\'{o}ri and Sz\'{e}kely (1983),
Andel and Dupas (1989), Martikainen and Petrov (1990), Petrov (2002), Xie (2008),
Feng, Li and Shen (2009), Frolov (2012) and references therein.
Note that M\'{o}ri and Sz\'{e}kely (1983) obtained the lower bound for
$\mathbf{P}(A_n \: i.o.)$ in terms of non-integer moments of
$\xi_n/\mathbf{E} \xi_n$. Generalizations of the first part of
the Borel--Cantelli lemma may be found in
Frolov (2014) and references therein.
One can find results under additional assumptions on dependence of events,
conditional Borel--Cantelli lemma and further references in
Chandra (2012).

The rest of the paper is organized as follows. We present our method and general
inequalities in section 2. Section 3 contains some new bounds
for probabilities of unions. Note that these results may be also applied
to measures of unions in arbitrary measurable spaces.
In section 4, new generalizations of the Borel--Cantelli lemma are proved.

\section{ Method and general results.}

Our method is based on the following result which is a generalization
of Theorem 1 in Frolov (2012).

\begin{theorem}\label{th1}
{\it
Let $\ell$ and $N$ be natural numbers such that $2 \leqslant \ell \leqslant N$.
Let $\{r_i, \, 1 \leqslant i \leqslant N\}$
and
$\{f_{ki},\, 1\leqslant k \leqslant \ell,\, 1 \leqslant i \leqslant N\}$
be arrays of non-negative real numbers. For $1\leqslant k \leqslant \ell$, put
\begin{eqnarray} \label{10}
\bar{s}_k = \sum\limits_{i=1}^N f_{ki} r_i.
\end{eqnarray}

Assume that there exist real numbers $c_1, c_2, \dots, c_N$ and
$a_1, a_2, \dots, a_\ell$ such that
\begin{eqnarray} \label{20}
\sum\limits_{i=1}^N (1-c_i) r_i = \sum\limits_{i=1}^\ell a_i \bar{s}_i.
\end{eqnarray}

If $c_i \geqslant 0$ for all $i=1,2,\ldots, N$, then
\begin{eqnarray} \label{30}
R = \sum\limits_{i=1}^N r_i \geqslant \sum\limits_{i=1}^\ell a_i \bar{s}_i.
\end{eqnarray}
If $c_i \leqslant 0$ for all $i=1,2,\ldots, N$, then
\begin{eqnarray} \label{40}
R \leqslant \sum\limits_{i=1}^\ell a_i \bar{s}_i.
\end{eqnarray}

Inequalities (\ref{30}) and (\ref{40}) turn to equalities if, for some
$1 \leqslant i_1<i_2<\dots < i_\ell \leqslant N$,
the numbers
$r_{i_1}, r_{i_2},\dots, r_{i_\ell}$ are solutions
of the linear system
\begin{eqnarray} \label{50}
&&
\sum\limits_{j=1}^\ell f_{k i_j} r_{i_j}=\bar{s}_k, \quad
1\leqslant k \leqslant \ell,
\end{eqnarray}
$r_{i}=0$ for all $i \neq i_k$
and $c_{i_k}=0$ for all $1\leqslant k \leqslant \ell$.
In this case, $R = r_{i_1}+ r_{i_2} +\cdots+ r_{i_\ell}$.
}
\end{theorem}

{\bf Proof.} By (\ref{20}), we have
$$ R - \sum\limits_{i=1}^\ell a_i \bar{s}_i = \sum\limits_{i=1}^N c_i r_i,
$$
that yields the assertion of Theorem \ref{th1}.
$\Box$

Now we describe the method which gives sharp lower and upper bounds for $R$.

We first choose number $\ell$ and array $\{f_{ki}\}$. The next step is to
take $c_i$. To satisfy relation (\ref{20}),
the simplest choice of $c_i$ is
$$ c_i = 1 - \sum\limits_{j=1}^\ell a_j f_{ji}
$$
for all $1\leqslant i \leqslant N$. Since the bounds for $R$ have to be sharp
(i.e. they have to turn to equalities for some set of numbers
$r_1, r_2, \ldots, r_N$),
we will take coefficients $a_j$ such that $c_{i_k}=0$ for
some $1 \leqslant i_1<i_2<\dots < i_\ell \leqslant N$.
To this end, we need a way to determine $i_1, i_2, \ldots, i_\ell$.

Note that if $f_{ki} = i^k$, $1 \leqslant k \leqslant \ell$, then
putting
$$ c_i = \prod\limits_{k=1}^\ell \left( 1- \frac{i}{i_k} \right)
$$
for all $1\leqslant i \leqslant N$ gives a simple way to find
$i_1, i_2, \ldots, i_\ell$ such that
$ c_i \geqslant 0$ or $ c_i \leqslant 0$ for every $i$.
Indeed, assume first that $\ell=2$. If we put $i_1 = m-1$ and $i_2=m$,
where $2\leqslant m \leqslant N$, then $ c_i \geqslant 0$ for all $i$.
If we take $i_1 = 1$ and $i_2 = N$, then $ c_i \leqslant 0$ for every $i$.
Here and in the sequel, natural number $m$ is a parameter which will be
specified in proofs below.
Further, suppose that $\ell=3$. If we choose $i_1 = m-1$, $i_2=m$ and $i_3=N$
with $2\leqslant m \leqslant N-1$, then $ c_i \geqslant 0$ for all $i$.
If we take $i_1 = 1$, $i_2=m-1$  and $i_3 = m$, where $3\leqslant m \leqslant N$,
then $ c_i \leqslant 0$ for every $i$. For $\ell=4$, putting
$i_1 = m-1$, $i_2=m$, $i_3=N-1$ and $i_4=N$ yields that
$ c_i \geqslant 0$ for all $i$. If $\ell=4$ and
$i_1 = 1$, $i_2 = m-1$, $i_3=m$  and $i_4 = N$, $4\leqslant m \leqslant N$,
then $ c_i \leqslant 0$ for every $i$.
And so on.
%
We will only use this way in the sequel to choose $i_1, i_2, \ldots, i_\ell$.

When we know $i_1, i_2, \ldots, i_\ell$, coefficients $a_j$ may be found
as solutions of the following linear system:
$$ \sum\limits_{j=1}^\ell a_j f_{ji_k} = 1, \quad k=1,2,\ldots,\ell.
$$

Taking into account that $f_{ki}$ may differ from $i^k$, we have to make certain
that the choice of $i_1, i_2, \ldots, i_\ell$ yields desired inequalities
for $c_i$. If we construct a lower bound, then we have to check
that $c_i \geqslant 0$ for all $i$. If we deal with an upper bound, then
we have to verify that $c_i \leqslant 0$ for all $i$.
By Theorem \ref{th1}, we get either inequality (\ref{30}), or
inequality (\ref{40}). Since indices $i_k$ depend on $m$,
we make an optimization over $m$.

Note that in the case $f_{ki} = i^k$, we get
$$ c_i = \prod\limits_{k=1}^\ell \left( 1- \frac{i}{i_k} \right)
= \sum \limits_{k=1}^\ell a_k i^k
$$
for all $i$. It means that $a_j$ are coefficients in the decomposition of
$c_i$ over degrees of $i$ and all $i, i^2, \ldots, i^\ell$ are in
this decomposition. Unfortunately, we can not follow this pattern to
find $a_j$ in general case. Indeed, if $f_{ki} = i^{\gamma_k}$,
$\gamma_k>0$, for example,
we could put
$$ c_i = \prod\limits_{k=1}^\ell \left( 1- \left(\frac{i}{i_k}
\right)^{\gamma_k-\gamma_{k-1}} \right).
$$
In this case we will not obtain a desired decomposition
with all $i^{\gamma_1}, i^{\gamma_2}, \ldots, i^{\gamma_\ell}$.
A variant with $\gamma_k$ instead of $\gamma_k-\gamma_{k-1}$
yields the same problem as well.

The above reasons lead us to the following result.

\begin{corollary}\label{c0}
Assume that for some $1 \leqslant i_1< i_2< \ldots < i_\ell \leqslant N$,
coefficients $a_j$ are solutions of the following linear system:
\begin{eqnarray}\label{51}
\sum\limits_{j=1}^\ell a_j f_{ji_k} = 1, \quad k=1,2,\ldots,\ell.
\end{eqnarray}
Put
\begin{eqnarray}\label{52}
c_i = 1 - \sum\limits_{j=1}^\ell a_j f_{ji}
\end{eqnarray}
for all $1\leqslant i \leqslant N$.
Assume that the numbers
$r_{i_1}^\ast, r_{i_2}^\ast,\dots, r_{i_\ell}^\ast$ are solutions
of the linear system (\ref{50}) and $r_{i}^\ast=0$ for all $i \neq i_k$,
$1 \leqslant i \leqslant N$.

If $c_i \geqslant 0$ for all $i=1,2,\ldots, N$, then
$R \geqslant R^\ast = r_{i_1}^\ast+ r_{i_2}^\ast +\cdots+ r_{i_\ell}^\ast$.

If $c_i \geqslant 0$ for all $i=1,2,\ldots, N$, then
$R \geqslant R^\ast$.
\end{corollary}

{\bf Proof.} By (\ref{52}), we have (\ref{20}).
If $c_i \geqslant 0$ for all $i=1,2,\ldots, N$, then inequality (\ref{30}) holds.
By definition of $r_{1}^\ast, r_{2}^\ast,\dots, r_{N}^\ast$, we get
\begin{eqnarray*}
\bar{s}_k^\ast = \sum\limits_{i=1}^N f_{ki} r_i^\ast = \bar{s}_k,
\quad k=1,2,\ldots,\ell.
\end{eqnarray*}
The last equality follows from (\ref{50}).
It yields that
\begin{eqnarray*}
R \geqslant \sum\limits_{i=1}^\ell a_i \bar{s}_i^\ast.
\end{eqnarray*}
By (\ref{51}), we have $c_{i_k}=0$ for all $k=1,2,\ldots,\ell$.
Applying Theorem \ref{th1} to $r_{1}^\ast, r_{2}^\ast,\dots, r_{N}^\ast$,
we conclude that
\begin{eqnarray*}
R^\ast = \sum\limits_{i=1}^\ell a_i \bar{s}_i^\ast \leqslant R.
\end{eqnarray*}
The case $c_i \leqslant 0$ for all $i=1,2,\ldots, N$ may be
considered in the same way.
$\Box$

We choose $\{f_{ik}\}$ such that results will be more simple.
To this end, in the sequel, we put $f_{ik} = i^{a+(k-1)\varrho}$ for all
$1\leqslant i \leqslant N$ and $1 \leqslant k \leqslant\ell$,
where $a>0$ and $\varrho >0$. Then relation (\ref{10}) turns to
\begin{eqnarray}\label{bsk}
\bar{s}_k = \sum\limits_{i=1}^N i^{a+(k-1)\varrho} r_i,
\quad 1 \leqslant k \leqslant\ell.
\end{eqnarray}
For the case $a=\varrho=1$, we will use a special notation
\begin{eqnarray}\label{sk}
s_k = \sum\limits_{i=1}^N i^{k} r_i,
\quad 1 \leqslant k \leqslant\ell.
\end{eqnarray}

We start with the case $\ell=2$. Our first result is a lower
bound for $R$.

\begin{theorem}\label{th2}
Define $\bar{s}_1$ and $\bar{s}_2$ by $(\ref{bsk})$.
Put $\bar{\delta}= (\bar{s}_2/\bar{s}_1)^{1/\varrho}$,
$\theta= \bar{\delta} -[\bar{\delta}]
$ and $\bar{\theta}=
(\bar{\delta}^{\rho}-(\bar{\delta}-\theta)^{\varrho})/
((\bar{\delta}+1-\theta)^{\varrho}-(\bar{\delta}-\theta)^{\varrho}) \in [0,1)
$,
where $[\cdot]$ denotes the integer part of the number in brackets.
Here and in the sequel, we suggest that $0/0=0$.

Then
\begin{eqnarray}\label{60}
R \geqslant
\frac{\bar{\theta} \bar{s}_1^{(a+\varrho)/\varrho} }{
\left(\bar{s}_2^{1/\varrho}+(1-\theta) \bar{s}_1^{1/\varrho}\right)^a}
+
\frac{(1-\bar{\theta}) \bar{s}_1^{(a+\varrho)/\varrho} }{
\left(\bar{s}_2^{1/\varrho}-\theta \bar{s}_1^{1/\varrho}\right)^a}.
\end{eqnarray}
\end{theorem}

Note that if $\bar{s}_1 = 0$, then $r_i=0$ for all $i$,
$\bar{s}_2 = \bar{\delta} = \theta = \bar{\theta} =0$ and
(\ref{60}) is trivial. If $\bar{s}_1>0$, then $\bar{\delta} \geqslant 1$.
Moreover, if $\bar{\delta} = 1$, then
$r_i=0$ for all $i\geqslant 2$, $\bar{s}_2 = \bar{s}_1 = r_1$,
$\theta = \bar{\theta} =0$  and (\ref{60}) turns to
$R \geqslant r_1$.

{\bf Proof.}
For every natural $m$, $2 \leqslant m\leqslant N$, put $i_1=m-1$ and $i_2=m$.
By (\ref{51}) and (\ref{52}), we have
$ c_i = 1 - a_1 i^a - a_2 i^{a+\varrho}$ for $i=1,2,\ldots, N$,
where $a_1$ and $a_2$ satisfy the following linear system:
\begin{eqnarray*}
&&
(m-1)^a a_1 + (m-1)^{a+\varrho} a_2 = 1, 
\\ &&
m^a a_1 + m^{a+\varrho} a_2 = 1.
\end{eqnarray*}
It follows that
$$ a_1 = \frac{m^{a+\varrho}-(m-1)^{a+\varrho}}{
m^{a+\varrho} (m-1)^a-m^a (m-1)^{a+\varrho}}, \quad
a_2 = -\frac{m^a-(m-1)^a}{m^{a+\varrho} (m-1)^a-m^a (m-1)^{a+\varrho}},
$$
and
$$ c_i = 1-\frac{m^{a+\varrho}-(m-1)^{a+\varrho}}{
m^{a+\varrho} (m-1)^a-m^a (m-1)^{a+\varrho}} i^a
+ \frac{m^a-(m-1)^a}{m^{a+\varrho} (m-1)^a-m^a (m-1)^{a+\varrho}} i^{a+\varrho}
$$
for all $i=1,2,\ldots, N$.

Now we check that $ c_i \geqslant 0$ for all $i$.
Consider function $f(x) = 1- a_1 x^a -a_2 x^{a+\varrho}$ for real $x \geqslant 0$.
We have $f(0)=1$, $f(m-1)=f(m)=0$, $f(+\infty) = +\infty$ and
$f'(x) = -  x^{a-1} (a a_1+ (a+\varrho) a_2 x^{\varrho})$. It is clear that
there exists a unique solution of equation $f'(x) = 0$. It
follows that function $f(x)$ takes its minimum at $x_0 \in (m-1, m)$
and $f(x) \leqslant 0$ for $x \in (m-1, m)$ and $f(x) \geqslant 0$ otherwise.
Hence, $ c_i = f(i) \geqslant 0$ for all $i=1,2,\ldots, N$.

Linear system (\ref{50}) is
\begin{eqnarray*}
&&
(m-1)^a r_{m-1}+ m^a r_{m} = \bar{s}_1, 
\\ &&
(m-1)^{a+\varrho} r_{m-1}+ m^{a+\varrho} r_{m} = \bar{s}_2,
\end{eqnarray*}
Solving this system, we get
\begin{eqnarray*}
r_{m-1}^\ast=\frac{\bar{s}_1 m^{a+\varrho} - \bar{s}_2 m^a }{
m^{a+\varrho} (m-1)^a-m^a (m-1)^{a+\varrho}},\quad
r_{m}^\ast=\frac{\bar{s}_2 (m-1)^a - \bar{s}_1 (m-1)^{a+\varrho} }{
m^{a+\varrho} (m-1)^a-m^a (m-1)^{a+\varrho}}.
\end{eqnarray*}

Inequalities $r_{m-1}^\ast\geqslant 0$ and $r_{m}^\ast\geqslant 0$ imply that
$(\bar{s}_2/\bar{s}_1)^{1/\varrho} \leqslant m
\leqslant 1+(\bar{s}_2/\bar{s}_1)^{1/\varrho}$,
which coincides with
$\bar{\delta} \leqslant m \leqslant 1+\bar{\delta}$.
This inequality defines $m$ uniquely for non-integer $\bar{\delta}$.
If $\bar{\delta}$ is integer then there are two variants.
Anyone of them maybe excluded. Hence,
without loss of generality, we assume in the sequel
that $\bar{\delta} < m \leqslant 1+\bar{\delta}$.
Then
$$ m = 1+ [\bar{\delta}] = 1+ \bar{\delta} - \theta,
$$
provided $\bar{\delta}<N$. It follows in this case that
$$ \bar{\theta} =
\frac{\bar{\delta}^\varrho - (m-1)^\varrho}{m^\varrho - (m-1)^\varrho} \in [0,1).
$$

Taking into account that $\bar{s}_2\leqslant N^\varrho \bar{s}_1$,
we see that $\bar{\delta}\leqslant N$. Hence,
we finally put $m=\min\{1+[\bar{\delta}],N\} \leqslant N$.

Assume that $\bar{\delta}<N$.  Then $m=1+[\bar{\delta}]$ and
$$ (m-1)^a = \bar{s}_1^{- a/\varrho}
\left(\bar{s}_2^{1/\varrho}-\theta \bar{s}_1^{1/\varrho} \right)^a .
$$
It follows that
$$ r_{m-1}^\ast = \frac{\bar{s}_1 ( m^{\varrho}-\bar{\delta}^{\varrho})}{
(m-1)^a (m^{\varrho} - (m-1)^{\varrho})}
= \frac{(1-\bar{\theta}) \bar{s}_1^{(a+\varrho)/\varrho} }{
\left(\bar{s}_2^{1/\varrho}-\theta \bar{s}_1^{1/\varrho}\right)^a}.
$$
Similarly,
$$ m^a = \bar{s}_1^{- a/\varrho}
\left(\bar{s}_2^{1/\varrho}+(1-\theta) \bar{s}_1^{1/\varrho} \right)^a .
$$
This implies that
$$ r_m^\ast = \frac{\bar{s}_1 ( \bar{\delta}^{\varrho} - (m-1)^{\varrho})}{
m^a (m^{\rho} - (m-1)^{\rho})} =
\frac{\bar{\theta} \bar{s}_1^{(a+\varrho)/\varrho} }{
\left(\bar{s}_2^{1/\varrho}+(1-\theta) \bar{s}_1^{1/\varrho}\right)^a}.
$$

By Corollary \ref{c0}, we have $R \geqslant r_{m-1}^\ast+r_m^\ast$.
Substituting of $r_{m-1}^\ast$ and $r_m^\ast$ in the last inequality
yields inequality (\ref{60}).

Assume that $\bar{\delta}=N$. Then
$$ 0 = \bar{s}_2 - N^{\varrho} \bar{s}_1
= \sum\limits_{i=1}^{N-1} (i^{\varrho} - N^{\varrho})  i^a r_i,
$$
where $(i^{\varrho} - N^{\varrho}) i^a < 0$ for $i \leqslant N-1$.
The latter implies that $r_1=r_2=\cdots = r_{N-1}=0$.
Moreover, $\bar{\delta}=N$ yields  that $\theta = \bar{\theta}=0$.
Hence inequality (\ref{60}) turns to inequality
$R \geqslant r_{N}$ which holds obviously.
$\Box$

Theorem \ref{th2} yields more simple bounds as follows.

\begin{corollary}\label{c1}
{\it Define $\bar{s}_1$ and $\bar{s}_2$ by $(\ref{bsk})$.
If $\varrho \geqslant 1$, then
\begin{eqnarray}\label{80}
R \geqslant
\frac{ \bar{s}_1^{(a+\varrho)/\varrho} }{\bar{s}_2^{a/\varrho}}.
\end{eqnarray}
If $\varrho < 1$, then
\begin{eqnarray}\label{81}
R \geqslant \frac{1-\bar{\theta}}{1-\theta}\,
\frac{ \bar{s}_1^{(a+\varrho)/\varrho} }{\bar{s}_2^{a/\varrho}}.
\end{eqnarray}
}
\end{corollary}

{\bf Proof.} Put
$$ g(\theta, \bar{\theta}) =
\frac{\bar{\theta} \bar{s}_1^{(a+\varrho)/\varrho}}{
\left(\bar{s}_2^{1/\varrho}+(1-\theta) \bar{s}_1^{1/\varrho}\right)^a}
+
\frac{(1-\bar{\theta}) \bar{s}_1^{(a+\varrho)/\varrho}}{
\left(\bar{s}_2^{1/\varrho}-\theta \bar{s}_1^{1/\varrho}\right)^a}.
$$
We have
$$ g(\theta, \bar{\theta}) = g(\theta, \theta) +
(\bar{\theta}-\theta)\left(
\frac{\bar{s}_1^{(a+\varrho)/\varrho} }{
\left(\bar{s}_2^{1/\varrho}+(1-\theta) \bar{s}_1^{1/\varrho}\right)^a}
-
\frac{ \bar{s}_1^{(a+\varrho)/\varrho} }{
\left(\bar{s}_2^{1/\varrho}-\theta \bar{s}_1^{1/\varrho}\right)^a}
\right).
$$
If $\bar{\theta} \leqslant \theta$, then
$g(\theta, \bar{\theta}) \geqslant g(\theta,\theta).$
It follows from inequality (\ref{60}) that
$$ R \geqslant
\frac{\theta \bar{s}_1^{(a+\varrho)/\varrho} }{
\left(\bar{s}_2^{1/\varrho}+(1-\theta) \bar{s}_1^{1/\varrho}\right)^a}
+
\frac{(1-\theta) \bar{s}_1^{(a+\varrho)/\varrho} }{
\left(\bar{s}_2^{1/\varrho}-\theta \bar{s}_1^{1/\varrho}\right)^a}.
$$
The right-hand side of the last inequality takes its minimum over
$\theta$ for $\theta=0$. Hence, inequality (\ref{60}) implies
inequality (\ref{80}) provided $\bar{\theta} \leqslant \theta$.

Inequality $\bar{\theta} \leqslant \theta$ is equivalent to
$\bar{\delta}^{\varrho} \leqslant (1-\theta) (\bar{\delta}-\theta)^{\varrho}
+ \theta (\bar{\delta}-\theta+1)^{\varrho}$. The latter inequality
holds when $x^{\varrho}$ is a convex function of $x$, i.e. for $\varrho\geqslant 1$.
For $\varrho < 1$, $x^{\varrho}$ is a concave function of $x$ and, therefore,
$\bar{\theta} > \theta$.

So, we have proved inequality (\ref{80}). To check inequality (\ref{81}),
we note that for $\bar{\theta} > \theta$,
$$ g(\theta, \bar{\theta}) \geqslant
\min\left\{\frac{\bar{\theta}}{\theta}, \frac{1-\bar{\theta}}{1-\theta}\right\}
g(\theta, \theta) = \frac{1-\bar{\theta}}{1-\theta} g(\theta, \theta).
$$
Inequality (\ref{81}) now follows from (\ref{60}).
$\Box$

Theorem \ref{th2} and Corollary \ref{c1} yield the following result in the case
$a=\varrho=1$.

\begin{corollary}\label{c2}
{\it Define $s_1$ and $s_2$ by (\ref{sk}).
Put $\delta = s_2/s_1$ and $\theta=\delta -[\delta]$.

The following inequality holds:
\begin{eqnarray*}
R \geqslant
\frac{\theta s_1^{2} }{s_2+(1-\theta) s_1} +
\frac{(1-\theta) s_1^{2} }{s_2-\theta s_1}
\geqslant \frac{s_1^{2} }{s_2}.
\end{eqnarray*}
}
\end{corollary}

Note that the middle part of the last inequality takes its minimum over
$\theta$ for $\theta=0$.

Now we turn to upper bounds for $\ell=2$. Our next result
is as follows.

\begin{theorem}\label{th3}
Define $\bar{s}_1$ and $\bar{s}_2$ by (\ref{bsk}).
The following inequality holds:
\begin{eqnarray}\label{110}
R \leqslant \frac{N^{a+\varrho}-1}{N^{a+\varrho}-N^a}  \bar{s}_1 -
\frac{N^a-1}{N^{a+\varrho}-N^a}  \bar{s}_2.
\end{eqnarray}
\end{theorem}

{\bf Proof.}
Take $i_1=1$ and $i_2=N$. By (\ref{51}) and (\ref{52}), we have
$c_i = 1 - a_1 i^a - a_2 i^{a+\varrho}$ for all $i=1,2,\ldots, N$,
where $a_1$ and $a_2$ are such that
\begin{eqnarray*}
&&
a_1 +  a_2 = 1, 
\\ &&
N^a a_1 + N^{a+\varrho} a_2 = 1.
\end{eqnarray*}
Hence
$$ a_1 = \frac{N^{a+\varrho}-1}{N^{a+\varrho} - N^a}, \quad
a_2 = -\frac{N^a-1}{N^{a+\varrho} - N^a},
$$
and
$$ c_i = 1- \frac{N^{a+\varrho}-1}{N^{a+\varrho} - N^a} i^a +
\frac{N^a-1}{N^{a+\varrho} - N^a} i^{a+\varrho}
$$
for all $i=1,2,\ldots, N$.

Let us check that $ c_i \leqslant 0$ for all $i$.
Consider again function $f(x) = 1- a_1 x^a -a_2 x^{a+\varrho}$
for real $x \geqslant 0$.
We have $f(0)=1$, $f(1)=f(N)=0$, $f(+\infty) = +\infty$ and
$f'(x) = -  x^{a-1} (a a_1+ (a+\varrho) a_2 x^{\varrho})$. We see that
there exists a unique solution $x_0$ of equation $f'(x) = 0$. It
follows that function $f(x)$ takes its minimum at $x_0 \in (1, N)$,
$f(x) \leqslant 0$ for $x \in (1, N)$ and $f(x) \geqslant 0$ otherwise.
Hence, $ c_i = f(i) \leqslant 0$ for all $i=1,2,\ldots, N$.

Therefore Theorem \ref{th1} yields inequality (\ref{110}).
$\Box$

Theorem \ref{th3} implies the following result for
$a=\varrho=1$.

\begin{corollary}\label{c4}
{\it Define $s_1$ and $s_2$ by (\ref{sk}).
The following inequality holds:
\begin{eqnarray*}
R \leqslant
\frac{N+1}{N} s_1 - \frac{1}{N} s_2.
\end{eqnarray*}
}
\end{corollary}

Comparing Theorems \ref{th2} and \ref{th3}, we see that
lower bounds seem more interesting for $\ell=2$.
The situation will change in the case $\ell=3$, to which we turn now.
We start again with a lower bound for $R$.

\begin{theorem}\label{th4}
Define $\bar{s}_1$, $\bar{s}_2$ and $\bar{s}_3$ by (\ref{bsk}).
Put $\bar{\delta}_1=N^{\varrho}\bar{s}_1-\bar{s}_2$,
$\bar{\delta}_2=N^{\varrho}\bar{s}_2-\bar{s}_3$,
$\bar{\delta}= (\bar{\delta}_2/\bar{\delta}_1)^{1/\varrho}$,
$\theta= \bar{\delta} -[\bar{\delta}]
$ and $\bar{\theta}=
(\bar{\delta}^{\varrho}-(\bar{\delta}-\theta)^{\varrho})/
((\bar{\delta}+1-\theta)^{\varrho}-(\bar{\delta}-\theta)^{\varrho}) \in [0,1)
$.

The following inequality holds:
\begin{eqnarray} \label{130}
R \geqslant
\frac{  \bar{\delta}_1 (1- \bar{\theta}) (N^a-(\bar{\delta}-\theta)^a)
}{N^a (\bar{\delta}-\theta)^a (N^\varrho-(\bar{\delta}-\theta)^\varrho)
}
+
\frac{  \bar{\delta}_1 \bar{\theta} (N^a-(\bar{\delta}-\theta+1)^a)
}{N^a (\bar{\delta}-\theta+1)^a (N^\varrho-(\bar{\delta}-\theta+1)^\varrho)
}
+\frac{\bar{s}_1}{N^a}.
\end{eqnarray}
\end{theorem}

Note that if $\bar{\delta}_1 =0$, then $r_1=r_2=\cdots = r_{N-1}=0$,
$\bar{\delta}=\bar{\delta}_2=\theta=\bar{\theta}=0$
and (\ref{130}) turns to $R \geqslant r_N$ that holds obviously.
If $\bar{\delta}_1 >0$, then $\bar{\delta} \geqslant 1$ in view of
$$ \bar{\delta}_1 = \sum\limits_{i=1}^{N-1} (N^\varrho - i^\rho) i^a r_i
\leqslant \sum\limits_{i=1}^{N-1} (N^\varrho - i^\rho) i^{a+\varrho} r_i =
\bar{\delta}_2.
$$
The latter also yields that $\bar{\delta} = 1$ only if
$r_2=\cdots = r_{N-1}=0$. In the last case,
$\bar{\delta}_1 = \bar{\delta}_2 = (N^\varrho-1) r_1$,
$\theta=\bar{\theta}=0$ and (\ref{130}) turns to $R \geqslant r_1+r_N$.

It follows that we may assume that $\bar{\delta} > 1$ in the sequel.

{\bf Proof.}
Take natural $m$, $2 \leqslant m\leqslant N-1$,
and put $i_1=m-1$, $i_2=m$, $i_3=N$.
By (\ref{51}) and (\ref{52}), we get
$ c_i= 1- a_{1} i^a - a_{2} i^{a+\varrho} - a_{3} i^{a+2 \varrho}$,
where $a_{1}$, $a_{2}$ and $a_{3}$ are determined by
the following system of linear equations
\begin{eqnarray*}
&&
(m-1)^a a_1 +  (m-1)^{a+\varrho} a_2 + (m-1)^{a+2 \varrho} a_3  = 1, 
\\ &&
m^a a_1 +  m^{a+\varrho} a_2 + m^{a+2 \varrho} a_3  = 1, 
\\ &&
N^a a_1 + N^{a+\varrho} a_2 + N^{a+2 \varrho} a_3 = 1.
\end{eqnarray*}
Then we have
\begin{eqnarray*}
&&
\hspace*{-1.1\parindent}
a_1 \!=\!
\frac{ m^{a+\varrho} N^{a+\varrho} (N^{\varrho} \!-\! m^{\varrho}) \!-\!
(m\!-\!1)^{a+\varrho} N^{a+\varrho} (N^{\varrho} \!-\! (m\!-1\!)^{\varrho})
\!+\! m^{a+\varrho} (m\!-\!1)^{a+\varrho} (m^{\varrho} \!-\! (m\!-\!1)^{\varrho})
}{\Delta_m},
\\ &&
\hspace*{-1.1\parindent}
a_2\! =\! -\frac{m^{a} N^{a} (N^{2\varrho}-m^{2\varrho})
\!-\! (m-1)^{a} N^{a} (N^{2\varrho}-(m-1)^{2\varrho})
\!+\! m^{a} (m-1)^{a} (m^{2\varrho}-(m-1)^{2\varrho})
}{\Delta_m},
\\ &&
\hspace*{-1.1\parindent}
a_3 \!=\! \frac{ m^{a} N^{a} (N^{\varrho}-m^{\varrho})
- (m-1)^{a} N^{a} (N^{\varrho}-(m-1)^{\varrho})
+ m^{a} (m-1)^{a} (m^{\varrho}-(m-1)^{\varrho})
}{\Delta_m},
\end{eqnarray*}
where
$\Delta_m = N^a m^a (m-1)^a (m^{\varrho}-(m-1)^{\varrho})
(N^{\varrho}-(m-1)^{\varrho})(N^{\varrho}-m^{\varrho})$.
Considering function
$f(x) = 1- a_{1} x^a - a_{2} x^{a+\varrho} - a_{3} x^{a+2 \varrho}$,
one can check that $ c_i \geqslant 0$ for all $i=1,2,\ldots, N$.

Solving linear system (\ref{50}), that is
\begin{eqnarray*}
&&
(m-1)^a r_{m-1} + m^a  r_{m}  +  N^a r_{N} = \bar{s}_1, 
\\ &&
(m-1)^{a+\varrho} r_{m-1} +  m^{a+\varrho} r_{m} +  N^{a+\varrho} r_{N}
= \bar{s}_2, 
\\ &&
(m-1)^{a+2 \varrho} r_{m-1} +  m^{a+2 \varrho} r_{m} + N^{a+2 \varrho} r_{N} = \bar{s}_3,
\end{eqnarray*}
we get
\begin{eqnarray*}
&&
\hspace*{-\parindent}
r_{m-1}^\ast = m^{a} N^{a} (N^{\varrho}-m^{\varrho})
\frac{ m^{\varrho} N^{\varrho} \bar{s}_1 -
(N^{\varrho}+m^{\varrho}) \bar{s}_2
+ \bar{s}_3
}{\Delta_m},
\\ &&
\hspace*{-\parindent}
r_m^\ast = - (m-1)^{a} N^{a} (N^{\varrho}-(m-1)^{\varrho})
\frac{ (m-1)^{\varrho} N^{\varrho} \bar{s}_1
- (N^{\varrho}+(m-1)^{\varrho}) \bar{s}_2
+ \bar{s}_3
}{\Delta_m},
\\ &&
\hspace*{-\parindent}
r_N^\ast = (m-1)^{a} m^{a} (m^{\varrho}-(m-1)^{\varrho})
\frac{ (m-1)^{\varrho} m^{\varrho} \bar{s}_1
- (m^{\varrho}+(m-1)^{\varrho}) \bar{s}_2
+ \bar{s}_3
}{\Delta_m}.
\end{eqnarray*}

Making use of inequalities $r_{m-1}^\ast\geqslant 0$ and $r_{m}^\ast\geqslant 0$,
we conclude that
$$ \left(\frac{\bar{\delta}_2}{\bar{\delta}_1}\right)^{1/\varrho} \leqslant
m \leqslant 1+\left(\frac{\bar{\delta}_2}{\bar{\delta}_1}\right)^{1/\varrho}.
$$
The latter is equivalent to $\bar{\delta} \leqslant m \leqslant 1+ \bar{\delta}$.
By the same reason as in the proof of Theorem \ref{th2}, we assume that
$\bar{\delta} < m \leqslant 1+ \bar{\delta}$.
Taking into account that $\bar{\delta}_2 \leqslant (N-1)^\varrho \bar{\delta}_1$,
we obtain $\bar{\delta} \leqslant N-1$. Hence, we put
$m=\min\{1+[\bar{\delta}], N-1\}$.

It is not difficult to check that
\begin{eqnarray*}
&&
r_{m-1}^\ast = \frac{  \bar{\delta}_1 (1- \bar{\theta})
}{ (\bar{\delta}-\theta)^a (N^\varrho-(\bar{\delta}-\theta)^\varrho)
},
\\ &&
r_m^\ast = \frac{  \bar{\delta}_1 \bar{\theta}
}{ (\bar{\delta}-\theta+1)^a (N^\varrho-(\bar{\delta}-\theta+1)^\varrho)
},
\\ &&
r_N^\ast = \frac{\bar{s}_1}{N^a} -
\frac{  \bar{\delta}_1}{N^a} \left(
\frac{  \bar{\theta}}{(N^\varrho-(\bar{\delta}-\theta+1)^\varrho)}
+
\frac{  1-\bar{\theta}}{(N^\varrho-(\bar{\delta}-\theta)^\varrho)}
\right).
\end{eqnarray*}
Corollary \ref{c0} implies that
$R \geqslant r_{m-1}^\ast+r_{m}^\ast+r_N^\ast$.
Substituting the formulae for $r_{m-1}^\ast$, $r_{m}^\ast$ and $r_N^\ast$ in
the last inequality, we arrive at (\ref{130}).
$\Box$

Theorem \ref{th4} implies the next result.

\begin{corollary}\label{c5}
Under notations of Theorem \ref{th4}, if $a\leqslant \varrho$, then
\begin{eqnarray} \label{131}
R \geqslant
\frac{  \bar{\delta}_1 (1- \bar{\theta}) (N^a-\bar{\delta}^a)
}{N^a (\bar{\delta}-\theta)^a (N^\varrho-\bar{\delta}^\varrho)
}
+
\frac{  \bar{\delta}_1 \bar{\theta} (N^a-(\bar{\delta}+1)^a)
}{N^a (\bar{\delta}-\theta+1)^a (N^\varrho-(\bar{\delta}+1)^\varrho)
}
+\frac{\bar{s}_1}{N^a}.
\end{eqnarray}
If $a\geqslant \varrho$, then
\begin{eqnarray} \label{132}
R \geqslant
\frac{  \bar{\delta}_1 (1- \bar{\theta}) (N^a-(\bar{\delta}-1)^a)
}{N^a (\bar{\delta}-\theta)^a (N^\varrho-(\bar{\delta}-1)^\varrho)
}
+
\frac{  \bar{\delta}_1 \bar{\theta} (N^a-\bar{\delta}^a)
}{N^a (\bar{\delta}-\theta+1)^a (N^\varrho-\bar{\delta}^\varrho)
}
+\frac{\bar{s}_1}{N^a}.
\end{eqnarray}
\end{corollary}

{\bf Proof.} We need the following technical result.

\begin{lemma}\label{l1}
If either $0<u<v<1$, or $1<u<v$,
then $f(x)= \frac{1-u^x}{1-v^x}$, $x>0$, is a decreasing function.
\end{lemma}

We omit the proof of  Lemma \ref{l1}.

%

Put $u=(\bar{\delta}-\theta)/N$ and $v=\bar{\delta}/N$.
Since $\bar{\delta}_2 \leqslant (N-1)^\varrho \bar{\delta}_1$, we
have $v<1$. If $a<\varrho$, then by Lemma \ref{l1},
$$ \frac{1-u^a}{1-v^a} > \frac{1-u^\varrho}{1-v^\varrho},
$$
which is equivalent to
$$\frac{N^a-(\bar{\delta}-\theta)^a}{N^a-\bar{\delta}^a} >
\frac{N^\varrho-(\bar{\delta}-\theta)^\varrho}{
N^\varrho-\bar{\delta}^\varrho}.
$$
We will have an opposite inequality for $a>\varrho$.

It follows that
$$
\frac{N^a-(\bar{\delta}-\theta)^a}{N^\varrho-(\bar{\delta}-\theta)^\varrho}
$$
take its minimum over $\theta$ for $\theta=0$, if $a \leqslant \varrho$,
and for $\theta=1$, if $a \geqslant \varrho$.

Making use of Lemma \ref{l1}, one can check that the same holds true for
$$
\frac{N^a-(\bar{\delta}-\theta+1)^a}{N^\varrho-(\bar{\delta}-\theta+1)^\varrho}.
$$
Now Corollary \ref{c5} follows from the latter and Theorem \ref{th4}.
$\Box$

For $\varrho \geqslant 1$, inequalities (\ref{131}) and (\ref{132}) imply
more simple bounds.

\begin{corollary}\label{c6}
Assume that $\varrho \geqslant 1$.
Then $\bar{\theta} \leqslant \theta$ and one can replace
$\bar{\theta}$ by $ \theta$ in (\ref{131}) and (\ref{132}).
Moreover, if $a\leqslant \varrho$ in addition, then
\begin{eqnarray*} 
R \geqslant
\frac{  \bar{\delta}_1 (N^a-\bar{\delta}^a)
}{N^a \bar{\delta}^a (N^\varrho-\bar{\delta}^\varrho)}
+\frac{\bar{s}_1}{N^a}.
\end{eqnarray*}
If $a\geqslant \varrho$ in addition, then
\begin{eqnarray*} 
R \geqslant
\frac{  \bar{\delta}_1  (N^a-(\bar{\delta}-1)^a)
}{N^a \bar{\delta}^a (N^\varrho-(\bar{\delta}-1)^\varrho)
}
+\frac{\bar{s}_1}{N^a}.
\end{eqnarray*}
\end{corollary}

{\bf Proof.} One can check that $\bar{\theta} \leqslant \theta$ and
one can put $\bar{\theta} = \theta =0$ in (\ref{131}) and (\ref{132})
in the same way as in the proof of Corollary \ref{c1}.
$\Box$

The next result is an upper bound for $R$.

\begin{theorem}\label{th5}
Define $\bar{s}_1$, $\bar{s}_2$ and $\bar{s}_3$ by (\ref{bsk}).
Put $\hat{\delta}_1=\bar{s}_2-\bar{s}_1$,
$\hat{\delta}_2= \bar{s}_3-\bar{s}_2$,
$\hat{\delta} =(\hat{\delta}_2/\hat{\delta}_1)^{1/\varrho}$,
$\theta= \hat{\delta} -[\hat{\delta}]$ and
$\hat{\theta}=
(\hat{\delta}^{\varrho}-(\hat{\delta}-\theta)^{\varrho})/
((\hat{\delta}+1-\theta)^{\varrho}-(\hat{\delta}-\theta)^{\varrho}) \in [0,1)
$.

The following inequality holds:
\begin{eqnarray} \label{140}
R \leqslant
\bar{s}_1-
\frac{\hat{\delta}_1 (1-\hat{\theta}) ((\hat{\delta}-\theta)^a-1)}{
(\hat{\delta}-\theta)^a
((\hat{\delta}-\theta)^\varrho -1)}
-
\frac{\hat{\delta}_1 \hat{\theta} ((\hat{\delta}-\theta+1)^a-1)}{
 (\hat{\delta}-\theta+1)^a
((\hat{\delta}-\theta+1)^\varrho -1)}.
\end{eqnarray}
\end{theorem}

Note that if $\hat{\delta}_1=0$, then then $r_i =0$ for all $i\geqslant 2$,
$\hat{\delta}_2=\hat{\delta}=\theta=\hat{\theta}=0$, $\bar{s}_1 =r_1$
and (\ref{140}) holds.
If $\hat{\delta}_1>0$, then taking into account that
$$ \hat{\delta}_2=
s_3-s_2 = \sum\limits_{i=2}^N i^{a+\varrho} (i^\varrho -1) r_i
\geqslant 2^{\varrho} \sum\limits_{i=2}^N i^{a} (i^\varrho -1) r_i =
2^{\varrho} (s_2-s_1)= 2^{\varrho} \hat{\delta}_1,
$$
we arrive at $\hat{\delta}\geqslant 2$.

{\bf Proof.}
Take natural $m$, $3 \leqslant m\leqslant N$,
and put $i_1=1$, $i_2 = m-1$, $i_3=m$.
By (\ref{51}) and (\ref{52}),
$ c_i= 1- a_{1} i^a - a_{2} i^{a+\varrho} - a_{3} i^{a+2 \varrho},
$
where $a_{1}$, $a_{2}$ and $a_{3}$ are solutions of linear system
\begin{eqnarray*}
&&
 a_1 + a_2 + a_3 = 1,
\\ &&
(m-1)^a a_1 +  (m-1)^{a+\varrho} a_2 + (m-1)^{a+2 \varrho} a_3  = 1, 
\\ &&
m^a a_1 +  m^{a+\varrho} a_2 + m^{a+2 \varrho} a_3  = 1.
\end{eqnarray*}
Then we have
\begin{eqnarray*}
&&
\hspace*{-\parindent}
a_1 =
\frac{ (m^{a+2\varrho}-1)((m-1)^{a+\varrho}-1) -
(m^{a+\varrho}-1)((m-1)^{a+2\varrho}-1)
}{\Delta_m},
\\ &&
\hspace*{-\parindent}
a_2 = -\frac{(m^{a+2\varrho}-1)((m-1)^{a}-1) - (m^a-1)((m-1)^{a+2\varrho}-1)
}{\Delta_m},
\\ &&
\hspace*{-\parindent}
a_3 = \frac{(m^{a+\varrho}-1)((m-1)^{a}-1) - (m^a-1)((m-1)^{a+\varrho}-1)
}{\Delta_m},
\end{eqnarray*}
where
$\Delta_m =  m^a (m-1)^a (m^{\varrho}-(m-1)^{\varrho})
((m-1)^{\varrho}-1)(m^{\varrho}-1)$.
Considering function
$f(x) = 1- a_{1} x^a - a_{2} x^{a+\varrho} - a_{3} x^{a+2 \varrho}$,
one can check that $ c_i \geqslant 0$ for all $i=1,2,\ldots, N$.

Linear system (\ref{50}) is as follows.
\begin{eqnarray*}
&&
r_1+ (m-1)^a r_{m-1} + m^a  r_{m}  = \bar{s}_1, 
\\ &&
r_1+(m-1)^{a+\varrho} r_{m-1} +  m^{a+\varrho} r_{m} = \bar{s}_2, 
\\ &&
r_1+(m-1)^{a+2 \varrho} r_{m-1} +  m^{a+2 \varrho} r_{m} = \bar{s}_3,
\end{eqnarray*}
We have
\begin{eqnarray*}
&& \hspace*{-\parindent}
r_1^\ast = m^{a} (m-1)^{a} (m^{\varrho} - (m-1)^{\varrho})
\frac{ m^{\varrho} (m-1)^{\varrho} \bar{s}_1
- (m^{\varrho}+(m-1)^{\varrho}) \bar{s}_2 + \bar{s}_3
}{\Delta_m},
\\ &&
\hspace*{-\parindent}
r_{m-1}^\ast = - m^{a} (m^{\varrho} - 1)
\frac{ m^{\varrho} \bar{s}_1 - (m^{\varrho}+1) \bar{s}_2 + \bar{s}_3
}{\Delta_m},
\\ &&
\hspace*{-\parindent}
r_{m}^\ast = (m-1)^{a} ((m-1)^{\varrho}-1)
\frac{ (m-1)^{\varrho} \bar{s}_1 - ((m-1)^{\varrho} + 1) \bar{s}_2 + \bar{s}_3
}{\Delta_m}.
\end{eqnarray*}

Again making use of $r_{m-1}^\ast\geqslant 0$ and $r_{m}^\ast\geqslant 0$,
we get
$$ \left(\frac{\hat{\delta}_2}{\hat{\delta}_1}\right)^{1/\varrho} \leqslant
m \leqslant 1+\left(\frac{\hat{\delta}_2}{\hat{\delta}_1}\right)^{1/\varrho}.
$$
This inequality is equivalent to
$\hat{\delta} \leqslant m \leqslant \hat{\delta}+1$.
By the same reason as in the proof of Theorem \ref{th2}, we may assume
that $\hat{\delta} < m \leqslant \hat{\delta}+1$.
Remember that $\hat{\delta}\geqslant 2$.
It follows that we can put
$m=\min\{1+[\hat{\delta}], N\}$.

It is not difficult to check that
\begin{eqnarray*}
&&
r_1^\ast = \bar{s}_1 - \hat{\delta}_1
\left(\frac{1-\hat{\theta}}{(\hat{\delta}-\theta)^\varrho -1}+
\frac{\hat{\theta}}{(\hat{\delta}+1-\theta)^\varrho -1}
\right),
\\ &&
r_{m-1}^\ast = \frac{\hat{\delta}_1 (1-\hat{\theta})}{
(\hat{\delta}-\theta)^a ((\hat{\delta}-\theta)^\varrho -1)},
\\ &&
r_{m}^\ast = \frac{\hat{\delta}_1 \hat{\theta}}{
(\hat{\delta}+1-\theta)^a ((\hat{\delta}+1-\theta)^\varrho -1)}.
\end{eqnarray*}

It follows from Corollary \ref{c0} that
$ R \leqslant r_1^\ast + r_{m-1}^\ast + r_{m}^\ast$.
Substituting of $r_1^\ast$, $r_{m-1}^\ast$ and $r_{m}^\ast$ in the latter
inequality yields  (\ref{140}).
$\Box$

Theorem \ref{th5} gives simpler bounds as well.

\begin{corollary}\label{c7}
Under notations of Theorem \ref{th5}, if $a\leqslant \varrho$, then
\begin{eqnarray*} 
R \leqslant
\bar{s}_1-
\frac{\hat{\delta}_1 (1-\hat{\theta}) (\hat{\delta}^a-1)}{
(\hat{\delta}-\theta)^a (\hat{\delta}^\varrho -1)}
-
\frac{\hat{\delta}_1 \hat{\theta} ((\hat{\delta} +1)^a-1)}{
 (\hat{\delta}-\theta+1)^a
((\hat{\delta}+1)^\varrho -1)}.
\end{eqnarray*}
If $a\geqslant \varrho$, then
\begin{eqnarray*} 
R \leqslant
\bar{s}_1-
\frac{\hat{\delta}_1 (1-\hat{\theta}) ((\hat{\delta}-1)^a-1)}{
(\hat{\delta}-\theta)^a
((\hat{\delta}-1)^\varrho -1)}
-
\frac{\hat{\delta}_1 \hat{\theta} (\hat{\delta}^a-1)}{
 (\hat{\delta}-\theta+1)^a (\hat{\delta}^\varrho -1)}.
\end{eqnarray*}
\end{corollary}

For $\varrho \geqslant 1$, Corollary \ref{c7} implies the next result.

\begin{corollary}\label{c8}
Assume that $\varrho \geqslant 1$.
If $a\leqslant \varrho$, then
\begin{eqnarray*} 
R \leqslant
\bar{s}_1-
\frac{\hat{\delta}_1 (\hat{\delta}^a-1)}{
\hat{\delta}^a (\hat{\delta}^\varrho -1)}.
\end{eqnarray*}
If $a\geqslant \varrho$, then
\begin{eqnarray*} 
R \leqslant
\bar{s}_1-
\frac{\hat{\delta}_1 ((\hat{\delta}-1)^a-1)}{
\hat{\delta}^a ((\hat{\delta}-1)^\varrho -1)}.
\end{eqnarray*}
\end{corollary}

Proofs of Corollaries \ref{c7} and \ref{c8} follow the same pattern as
those of Corollaries \ref{c5} and \ref{c6}. We omit details.

\begin{remark}\label{r1}
The right-hand side of (\ref{140}) may be obtained from
the right-hand side of (\ref{130}) by a formal replacement
$N=1$, $\bar{\delta}_1 = - \hat{\delta}_1$, $\bar{\delta}_2 = - \hat{\delta}_2$,
$\bar{\delta} =  \hat{\delta}$ and $ \bar{\theta} = \hat{\theta}$.
\end{remark}

\begin{remark}\label{r2}
All inequalities of Theorems \ref{th2} -- \ref{th5} are sharp.
For each of these inequalities, there exists a set of numbers
$r_1, r_2, \ldots, r_N$ such that the inequality turns to equality.
\end{remark}

\section{Bounds for probabilities of unions.}

In this section, we discuss bounds for probabilities of unions of events
which follow from the results of section 2. Note that these bounds maybe
applied to measures of unions of sets in arbitrary measurable spaces.

Let $(\Omega,\mathcal{F},\mathbf{P})$ be a probability space. For
events $A_1, A_2,\dots,A_N$, put $U=\bigcup\limits_{i=1}^N A_i$.
Denote $B_i=\{\omega \in \Omega : \omega \: \mbox{belongs to exactly}\: i
\:\mbox{events of}\: A_1, A_2,\dots,A_N\}$ and
$p_i= \mathbf{P}(B_i)$, $i=0, 1, \dots, N$.
Then
$$ \mathbf{P} (U) = \sum\limits_{i=1}^N p_i.
$$

The simplest application of the above method is to put $r_i=p_i$ for
$i=1, 2, \dots, N$. Then the general results of the previous sections
yield Theorems 2--5 and Corollaries 1 and 2 in Frolov (2012) and
Theorem 3 and 4 and Corollaries 1--5 in Frolov (2014)
that are generalizations of earlier results.

Note that one may also consider more general events that maybe
represented by sums $p_{i_1}+p_{i_2}+\cdots+p_{i_M}$
where $i_1<i_2<\cdots < i_M$, $M \leqslant N$.
For example, sum $p_{t}+p_{t+1}+\cdots+p_{N}$
equals to probability that at least $t$ events from
$A_1, A_2,\dots,A_N$ occur. This requires a modification
of the above method and will be done elsewhere.

We now turn to another representations of $\mathbf{P} (U)$
which is a start point of our method as well.
By Lemma 1 in Kuai, Alajaji and Takahara (2000), we have
$$ \mathbf{P} (U) = \sum\limits_{k=1}^N \sum\limits_{i=1}^N
\frac{1}{i} p_{ik},
$$
where $p_{ik} = \mathbf{P} (B_i A_k)$.

We also give a simple proof of the last equality.
Putting $\xi_N=I_{A_1}+I_{A_2}+\cdots+I_{A_N}$,
we get $\xi_N I_{B_i} = i I_{B_i}$ and
$$ \sum\limits_{k=1}^N \sum\limits_{i=1}^N \frac{1}{i} p_{ik}
= \mathbf{E} \left( \sum\limits_{k=1}^N \sum\limits_{i=1}^N
\frac{1}{i} I_{B_i} I_{A_k} \right) =
\mathbf{E} \left( \sum\limits_{i=1}^N
\frac{1}{i} I_{B_i} \xi_N \right) =
\mathbf{E} \left( \sum\limits_{i=1}^N
I_{B_i} \right) = \mathbf{P} (U).
$$

For every fixed $k$, putting $r_{ik} = p_{ik}/i$ and
$$ R_k = \sum\limits_{i=1}^N r_{ik},
$$
we can take bounds for $R_k$ from our general results.

Denote
\begin{eqnarray}\label{bskk}
&&
\bar{s}_j(k) = \sum\limits_{i=1}^N i^{a+(j-1)\varrho} r_{ik},
\\ &&
\label{skk}
s_j(k) = \sum\limits_{i=1}^N i^{j} r_{ik},
\end{eqnarray}
for $1\leqslant j \leqslant \ell$ and $k=1,2,\ldots, N$.

An application of Theorem \ref{th2} yields the next result.

\begin{theorem}\label{th6}
Define $ \bar{s}_1(k)$ and $ \bar{s}_2(k)$ by (\ref{bskk}),
$k=1,2,\ldots, N$.
Put $\bar{\delta}_k= (\bar{s}_2(k)/\bar{s}_1(k))^{1/\varrho}$,
$\theta_k= \bar{\delta}_k -[\bar{\delta}_k]
$ and $\bar{\theta}_k=
(\bar{\delta}_k^{\varrho}-(\bar{\delta}_k-\theta_k)^{\varrho})/
((\bar{\delta}_k+1-\theta_k)^{\varrho}-(\bar{\delta}_k-\theta_k)^{\varrho}) \in [0,1)
$, where $k=1,2,\ldots, N$.

Then
\begin{eqnarray}\label{200}
\mathbf{P} (U) \geqslant \sum\limits_{k=1}^N
\left\{
\frac{\bar{\theta}_k \bar{s}_1^{(a+\varrho)/\varrho}(k) }{
\left(\bar{s}_2^{1/\varrho}(k)+(1-\theta_k) \bar{s}_1^{1/\varrho}(k)\right)^a}
+
\frac{(1-\bar{\theta}_k) \bar{s}_1^{(a+\varrho)/\varrho}(k) }{
\left(\bar{s}_2^{1/\varrho}(k)-\theta_k \bar{s}_1^{1/\varrho}(k)\right)^a}
\right\}.
\end{eqnarray}
\end{theorem}

For $a=\varrho=1$, Theorem \ref{th6} implies Theorem 1 in
Kuai, Alajaji and Takahara (2000).
By Corollary \ref{c2}, we may put $\bar{\theta}_k= \theta_k =0$
in (\ref{200}) and obtain a result in de Caen (1997).

It is clear that one can use all results from section 2 to
derive upper and lower bounds similar to that of Theorem \ref{th6}.
Theorem \ref{th4} implies the following result.

\begin{theorem}\label{th7}
Define $ \bar{s}_1(k)$, $ \bar{s}_2(k)$ and $ \bar{s}_3(k)$ by (\ref{bskk}),
$k=1,2,\ldots, N$.
Put $\bar{\delta}_{1k}=N^{\varrho}\bar{s}_1(k)-\bar{s}_2(k)$,
$\bar{\delta}_{2k}=N^{\varrho}\bar{s}_2(k)-\bar{s}_3(k)$,
$\bar{\delta}_k= (\bar{\delta}_{2k}/\bar{\delta}_{1k})^{1/\varrho}$,
$\theta_k= \bar{\delta}_k -[\bar{\delta}_k]
$ and $\bar{\theta}_k=
(\bar{\delta}_k^{\varrho}-(\bar{\delta}_k-\theta_k)^{\varrho})/
((\bar{\delta}_k+1-\theta_k)^{\varrho}-(\bar{\delta}_k-\theta_k)^{\varrho})
\in [0,1)$, $k=1,2,\ldots,N$.

The following inequality holds:
\begin{eqnarray*} 
\mathbf{P} (U) \!\geqslant\!\! \sum\limits_{k=1}^N \!
\left\{ \!
\frac{  \bar{\delta}_{1k} (1 \!-\! \bar{\theta}_k) (N^a \!- \!(\bar{\delta}_k\!-\!\theta_k)^a)
}{N^a (\bar{\delta}_k\!-\!\theta_k)^a (N^\varrho \!- \!(\bar{\delta}_k\!-\!\theta_k)^\varrho)
}
+
\frac{  \bar{\delta}_{1k} \bar{\theta}_k
(N^a \!- \!(\bar{\delta}_k\!-\!\theta_k \!+ \!1)^a)
}{N^a (\bar{\delta}_k\!-\!\theta_k \!+\! 1)^a
(N^\varrho-(\bar{\delta}_k\!-\!\theta_k\!+\! 1)^\varrho)
}
+\frac{\bar{s}_1(k)}{N^a}
\!\right\}\!.
\end{eqnarray*}
\end{theorem}

For $a=\varrho=1$, we obtain the next result from Corollary \ref{c6}.

\begin{corollary}\label{c10}
Define $ s_1(k)$, $ s_2(k)$ and $ s_3(k)$ by (\ref{skk})
and put
$\bar{\delta}_{1k}=N s_1(k)-s_2(k)$, $\bar{\delta}_{2k}=N s_2(k)- s_3(k)$
for $k=1,2,\ldots,N$.

The following inequality holds:
\begin{eqnarray} \label{211}
\mathbf{P} (U) \geqslant \frac{1}{N}
\sum\limits_{k=1}^N
\left\{
\frac{  \bar{\delta}_{1k}^2}{\bar{\delta}_{2k}} + s_1(k)
\right\}.
\end{eqnarray}
\end{corollary}

Note that for all $k=1,2,\ldots,N$, we have
\begin{eqnarray*}
&&
s_1(k) = \sum\limits_{i=1}^N p_{ik}
= \mathbf{E} \left( \sum\limits_{i=1}^N I_{B_i} I_{A_k} \right) =
\mathbf{E} \left( I_{U}  I_{A_k} \right) = \mathbf{P} (A_k),
\\ &&
s_2(k) = \sum\limits_{i=1}^N i p_{ik}
= \mathbf{E} \left( \sum\limits_{i=1}^N i I_{B_i} I_{A_k} \right) =
\mathbf{E} \left( \xi_N  I_{A_k} \right) = \sum\limits_{i=1}^N \mathbf{P} (A_i A_k),
\\ &&
s_3(k) = \sum\limits_{i=1}^N i^2 p_{ik}
= \mathbf{E} \left( \sum\limits_{i=1}^N i^2 I_{B_i} I_{A_k} \right) =
\mathbf{E} \left( \xi_N^2  I_{A_k} \right) =
\sum\limits_{i=1}^N \sum\limits_{j=1}^N \mathbf{P} (A_i A_j A_k).
\end{eqnarray*}
It follows that
\begin{eqnarray}\label{bdk}
\bar{\delta}_{1k} \!=\! \sum\limits_{i=1}^N \mathbf{P} (\overline{A_i} A_k)
\!=\! \mathbf{E} \left( (N \!-\! \xi_N)  I_{A_k} \right)
,\;
\bar{\delta}_{2k} \!=\! \sum\limits_{i=1}^N \sum\limits_{j=1}^N
\mathbf{P} (A_i \overline{A_j} A_k) \!=\!
\mathbf{E} \left( \xi_N (N \!-\! \xi_N) I_{A_k} \right),
\end{eqnarray}
for all $k=1,2,\ldots,N$.

Now we turn to upper bounds. The next result follows from Theorem \ref{th5}.

\begin{theorem}\label{th8}
Define $\bar{s}_1(k)$, $\bar{s}_2(k)$ and $\bar{s}_3(k)$ by (\ref{bskk}),
$k=1,2,\ldots, N$.
Put $\hat{\delta}_{1k}=\bar{s}_2(k)-\bar{s}_1(k)$,
$\hat{\delta}_{2k}= \bar{s}_3(k)-\bar{s}_2(k)$,
$\hat{\delta}_k =(\hat{\delta}_{2k}/\hat{\delta}_{1k})^{1/\varrho}$,
$\theta_k= \hat{\delta}_k -[\hat{\delta}_k]$ and
$\hat{\theta}_k=
(\hat{\delta}_k^{\varrho}-(\hat{\delta}_k-\theta_k)^{\varrho})/
((\hat{\delta}_k+1-\theta_k)^{\varrho}-(\hat{\delta}_k-\theta_k)^{\varrho}) \in [0,1)
$, $k=1,2,\ldots,N$.

The following inequality holds:
\begin{eqnarray*} 
\mathbf{P} (U) \leqslant \sum\limits_{k=1}^N
\left\{
\bar{s}_1(k)-
\frac{\hat{\delta}_{1k} (1-\hat{\theta}_k) ((\hat{\delta}_k-\theta_k)^a-1)}{
(\hat{\delta}_k-\theta_k)^a
((\hat{\delta}_k-\theta_k)^\varrho -1)}
-
\frac{\hat{\delta}_{1k} \hat{\theta}_k ((\hat{\delta}_k-\theta_k+1)^a-1)}{
 (\hat{\delta}_k-\theta_k+1)^a
((\hat{\delta}_k-\theta_k+1)^\varrho -1)}
\right\}.
\end{eqnarray*}
\end{theorem}

For $a=\varrho=1$, we obtain the next result from Corollary \ref{c8}.

\begin{corollary}\label{c11}
Define $ s_1(k)$, $ s_2(k)$ and $ s_3(k)$ by (\ref{skk})
and put $\hat{\delta}_{1k}= s_2(k)- s_1(k)$,
$\hat{\delta}_{2k}= s_3(k)-s_2(k)$
for $k=1,2,\ldots,N$.

The following inequality holds:
\begin{eqnarray} \label{221}
\mathbf{P} (U) \leqslant \sum\limits_{k=1}^N
\left\{
s_1(k)- \frac{\hat{\delta}_{1k}^2}{\hat{\delta}_{2k}}
\right\}.
\end{eqnarray}
\end{corollary}

Note that
\begin{eqnarray} \label{hdk1}
&&
\hat{\delta}_{1k} = \sum\limits_{i=1}^N \mathbf{P} (A_i A_k) - \mathbf{P} (A_k) =
\mathbf{E} (\xi_N - 1) I_{A_k}, 
\\ &&
\label{hdk2}
\hat{\delta}_{2k}
= \sum\limits_{i=1}^N \sum\limits_{j=1}^N \mathbf{P} (A_i A_j A_k) -
\sum\limits_{i=1}^N \mathbf{P} (A_i A_k) =
\mathbf{E} \xi_N (\xi_N - 1) I_{A_k}.
\end{eqnarray}

We finally mention that Theorems \ref{th6}--\ref{th8} and Corollaries
\ref{c10} and \ref{c11} are new result.

\section{Borel--Cantelli lemmas.}

Let $(\Omega,\mathcal{F},\mathbf{P})$ be a probability space and
$\{A_n\}$ be a sequence of events. Denote
$$ \{ A_n \;\mbox{i.o.}\} = \limsup A_n =
\bigcap\limits_{n=1}^\infty \bigcup\limits_{k=n}^\infty A_k.
$$
For $m \leqslant m$, put $U_{mn} = \bigcup\limits_{k=m}^n A_k$.
Since
$$ \mathbf{P}\left( A_n \;\mbox{i.o.}\right) =
\lim\limits_{m\to \infty}\lim\limits_{n\to \infty}
\mathbf{P}\left( U_{mn}
\right),
$$
every new upper or lower bound allows us to derive new variant
of first or second part of the Borel--Cantelli Lemma.
Our results of the previous section imply that
$$  Q(m,n) \leqslant \mathbf{P}\left( U_{mn}
\right)
\leqslant Q'(m,n),
$$
where $Q(m,n)$ and $Q'(m,n)$ are
the right-hand sides of the applied lower and upper bound,
correspondingly. It is clear that
$$
\limsup\limits_{m\to \infty}\limsup\limits_{n\to \infty} Q(m,n) \leqslant
\mathbf{P}\left( A_n \;\mbox{i.o.}\right) \leqslant
\liminf\limits_{m\to \infty}\liminf\limits_{n\to \infty} Q'(m,n).
$$
It may happen that we cannot find these double limits.
But, if for every fixed $m$ the inequality
$$ \limsup\limits_{n\to \infty} Q(m,n) \geqslant
\limsup\limits_{n\to \infty} Q(1,n)
$$
holds, then we have
$$ \mathbf{P}\left( A_n \;\mbox{i.o.}\right) \geqslant
\limsup\limits_{n\to \infty} Q(1,n).
$$
Similarly, if for every fixed $m$ the inequality
$$ \liminf\limits_{n\to \infty} Q'(m,n) \leqslant
\liminf\limits_{n\to \infty} Q'(1,n)
$$
holds, then we get
$$ \mathbf{P}\left( A_n \;\mbox{i.o.}\right) \leqslant
\liminf\limits_{n\to \infty} Q'(1,n).
$$

The most applicable variants of the Borel--Cantelli lemma are proved
by this way.
In the proof of the first part of the classical Borel--Cantelli lemma,
the upper bound for $ \mathbf{P}\left( U_{mn}\right)$ by
$s_1$ is used. In the Erd\H{o}s--R\'{e}nyi generalization
of the second part of the Borel--Cantelli lemma,
the inequality with $s_1$ and $s_2$ is applied.
Frolov (2012) has applied the lower bound
for $ \mathbf{P}\left( U_{mn}\right)$, based on
$s_1$, $s_2$ and $s_3$. This yielded a generalization
of the second part of the Borel--Cantelli lemma in
Theorem 9 of the last paper. Note that mentioned here bounds
are constructed for probabilities $r_i=p_i$.

In this section, we present new variants of the Borel--Cantelli lemma based on
inequalities (\ref{211}) and (\ref{221}). Note that
the last inequalities are constructed from bounds for numbers $r_{ik}=p_{ik}/i$.

We start with the second part of the Borel--Cantelli lemma.

\begin{theorem}\label{th9}
Denote $\xi_{n} = I_{A_1}+I_{A_2}+\cdots+I_{A_n}$
and $\eta_n= n-\xi_{n}$ for all natural $n$.

Assume that
$$ \frac{1}{n} \sum\limits_{k=1}^n
\frac{  \mathbf{E} \eta_n I_{A_k}}{
\mathbf{E} \eta_n \xi_n I_{A_k}} \rightarrow 0 \quad\mbox{as}\quad
n \rightarrow\infty.
$$

Then
$$ \mathbf{P}\left( A_n \;\mbox{i.o.}\right)
\geqslant \limsup\limits_{n\to \infty}
\frac{1}{n} \sum\limits_{k=1}^n \left\{
\mathbf{P}\left( A_k \right)+
\frac{ (\mathbf{E} \eta_n I_{A_k})^2}{\mathbf{E} \eta_n \xi_n I_{A_k}}
\right\}.
$$

\end{theorem}

It follows from (\ref{bdk}) that
$$ \mathbf{E} \eta_n I_{A_k} = \sum\limits_{i=1}^n
\mathbf{P}\left( \overline{A_i} A_k \right), \quad
\mathbf{E} \eta_n \xi_n I_{A_k} =
\sum\limits_{i=1}^n
\sum\limits_{j=1}^n \mathbf{P}\left( \overline{A_i} A_j A_k \right).
$$

{\bf Proof.}
Inequality (\ref{211}) and relation (\ref{bdk}) yield that
$$ \mathbf{P}\left( \bigcup\limits_{k=m}^n A_k \right) \geqslant
\frac{1}{n-m+1}\sum\limits_{k=m}^n
\left\{
\mathbf{P}\left( A_k \right)+ T_k(m,n) \right\},
$$
where
$$ T_k(m,n) = \frac{ \left( \mathbf{E} (\eta_n - \eta_{m-1}) I_{A_k} \right)^2}{
\mathbf{E} (\eta_n - \eta_{m-1})(\xi_n - \xi_{m-1}) I_{A_k}}.
$$

We have
\begin{eqnarray*}
&&
T_k(m,n) 
\geqslant
\frac{ \left( \mathbf{E} (\eta_n - \eta_{m-1}) I_{A_k} \right)^2}{
\mathbf{E} \eta_n \xi_n I_{A_k}}
= \frac{ (\mathbf{E} \eta_n I_{A_k})^2 - 2 \mathbf{E} \eta_n \eta_{m-1} I_{A_k}
+(\mathbf{E}\eta_{m-1} I_{A_k})^2}{
\mathbf{E} \eta_n \xi_n I_{A_k} 
}
\\ &&
\geqslant
\frac{ (\mathbf{E} \eta_n I_{A_k})^2 - 2(m-1) \mathbf{E} \eta_n I_{A_k}}{
\mathbf{E} \eta_n \xi_n I_{A_k}
}.
\end{eqnarray*}
By (\ref{211}), the inequality
$$
1 \geqslant \frac{1}{n} \sum\limits_{k=1}^n T_k(1,n) =
\frac{1}{n} \sum\limits_{k=1}^n \frac{ (\mathbf{E} \eta_n I_{A_k})^2}{
\mathbf{E} \eta_n \xi_n I_{A_k}}
$$
holds for all natural $n$. It implies that
$$
\sum\limits_{k=1}^m
\frac{ (\mathbf{E} \eta_n I_{A_k})^2}{\mathbf{E} \eta_n \xi_n I_{A_k}}
\leqslant m.
$$

Hence
\begin{eqnarray*}
&&
\hspace*{-\parindent}
\mathbf{P}\left( \bigcup\limits_{k=m}^n A_k \right) \geqslant
\frac{1}{n} \sum\limits_{k=m}^n
\left\{
\mathbf{P}\left( A_k \right)+
\frac{ (\mathbf{E} \eta_n I_{A_k})^2}{\mathbf{E} \eta_n \xi_n I_{A_k}} -
\frac{ 2(m-1) \mathbf{E} \eta_n I_{A_k}}{
\mathbf{E} \eta_n \xi_n I_{A_k}}\right\}
\\ && \hspace*{-\parindent}
\geqslant
\frac{1}{n} \sum\limits_{k=1}^n
\left\{
\mathbf{P}\left( A_k \right)+
\frac{ (\mathbf{E} \eta_n I_{A_k})^2}{\mathbf{E} \eta_n \xi_n I_{A_k}}
\right\} - \frac{2m}{n} -
\frac{2(m-1)}{n} \sum\limits_{k=1}^n
\frac{  \mathbf{E} \eta_n I_{A_k}}{
\mathbf{E} \eta_n \xi_n I_{A_k}}.
\end{eqnarray*}

This yields that for every fixed $m$,
$$ \mathbf{P}\left( \bigcup\limits_{k=m}^\infty A_k \right) \geqslant
\limsup\limits_{n \to \infty}
\frac{1}{n} \sum\limits_{k=1}^n
\left\{\mathbf{P}\left( A_k \right)+
\frac{ (\mathbf{E} \eta_n I_{A_k})^2}{\mathbf{E} \eta_n \xi_n I_{A_k}}
\right\}.
$$
The last inequality implies the desired assertion.
$\Box$

Theorem 9 in Frolov (2012) contains a lower bound for
$\mathbf{P}\left( A_n \;\mbox{i.o.}\right)$ constructed
from $p_i$. There is an example in Frolov (2012) which
shows that this lower bound is better than previous ones.
One can check that for this example, the lower bounds
of Theorem 9 in Frolov (2012) and Theorem \ref{th9} of this
section coincide.

Now we turn to the first part of the Borel--Cantelli lemma.

\begin{theorem}\label{th10}
Denote $\xi_{n} = I_{A_1}+I_{A_2}+\cdots+I_{A_n}$
for all natural $n$ and $\xi_0=0$.
If
\begin{eqnarray}\label{usl}
\sum\limits_{k=m}^n
\frac{  \mathbf{E} (\xi_n - \xi_{m-1}) I_{A_k}}{
\mathbf{E} (\xi_n- \xi_{m-1})^2 I_{A_k}} \rightarrow 0 \quad\mbox{as}\quad
n \rightarrow\infty
\end{eqnarray}
for all sufficiently large $m$, then
$$ \mathbf{P}\left( A_n \;\mbox{i.o.}\right)
\leqslant \limsup\limits_{m\to \infty} \limsup\limits_{n\to \infty}
\sum\limits_{k=m}^n \left\{
\mathbf{P}\left( A_k \right)-
\frac{ (\mathbf{E} (\xi_n- \xi_{m-1}) I_{A_k})^2}{\mathbf{E} (\xi_n- \xi_{m-1})^2 I_{A_k}}
\right\}.
$$


If condition (\ref{usl}) holds for $m=1$, then
$$ \mathbf{P}\left( A_n \;\mbox{i.o.}\right)
\leqslant \limsup\limits_{m\to \infty} \limsup\limits_{n\to \infty}
\sum\limits_{k=m}^n \left\{
\mathbf{P}\left( A_k \right)-
\frac{ (\mathbf{E} \xi_n I_{A_k})^2}{\mathbf{E} \xi_n^2 I_{A_k}}
\right\}.
$$

\end{theorem}

It follows from (\ref{hdk1}) and (\ref{hdk2}) that
$$ \mathbf{E} \xi_n I_{A_k} = \sum\limits_{i=1}^n
\mathbf{P}\left( A_i A_k \right), \quad
\mathbf{E}  \xi_n^2 I_{A_k} =
\sum\limits_{i=1}^n
\sum\limits_{j=1}^n \mathbf{P}\left(A_i A_j A_k \right).
$$

{\bf Proof.}
Inequality (\ref{221}) and assertions (\ref{hdk1}) and (\ref{hdk2}) imply that
$$ \mathbf{P}\left( \bigcup\limits_{k=m}^n A_k \right) \leqslant
\sum\limits_{k=m}^n
\left\{
\mathbf{P}\left( A_k \right)- T'_k(m,n) \right\},
$$
where
$$ T'_k(m,n) = \frac{ \left( \mathbf{E} (\xi_n - \xi_{m-1}-1) I_{A_k} \right)^2}{
\mathbf{E} (\xi_n - \xi_{m-1})(\xi_n - \xi_{m-1}-1) I_{A_k}}.
$$
We have
\begin{eqnarray*}
&&
T'_k(m,n) 
\geqslant
\frac{ \left( \mathbf{E} (\xi_n - \xi_{m-1}-1) I_{A_k} \right)^2}{
\mathbf{E} (\xi_n- \xi_{m-1})^2 I_{A_k}}
\geqslant
\frac{ \left( \mathbf{E} (\xi_n - \xi_{m-1}) I_{A_k} \right)^2
- 2 \mathbf{E} (\xi_n - \xi_{m-1}) I_{A_k} }{
\mathbf{E} (\xi_n- \xi_{m-1})^2 I_{A_k}}
\\ &&
\geqslant \frac{ (\mathbf{E} \xi_n I_{A_k})^2 -
2 \mathbf{E} \xi_n \xi_{m-1} I_{A_k} - 2 \mathbf{E}\xi_{n} I_{A_k}}{
\mathbf{E} \xi_n^2 I_{A_k} 
}
\geqslant
\frac{ (\mathbf{E} \xi_n I_{A_k})^2 - 2 m \mathbf{E} \xi_n I_{A_k}}{
\mathbf{E} \xi_n^2 I_{A_k} 
}.
\end{eqnarray*}
It yields that
\begin{eqnarray*}
&&
\hspace*{-\parindent}
\mathbf{P}\left( \bigcup\limits_{k=m}^n A_k \right) \leqslant
\sum\limits_{k=m}^n
\left\{
\mathbf{P}\left( A_k \right)-
\frac{ (\mathbf{E} (\xi_n- \xi_{m-1}) I_{A_k})^2}{
\mathbf{E} (\xi_n- \xi_{m-1})^2 I_{A_k}} \right\} +
2 \sum\limits_{k=m}^n
\frac{ \mathbf{E} (\xi_n- \xi_{m-1}) I_{A_k}}{
\mathbf{E}  (\xi_n- \xi_{m-1})^2 I_{A_k}}
\\ && \hspace*{-\parindent}
\leqslant\sum\limits_{k=m}^n
\left\{
\mathbf{P}\left( A_k \right)-
\frac{ (\mathbf{E} \xi_n I_{A_k})^2}{\mathbf{E} \xi_n^2 I_{A_k}}
\right\}+
2 m \sum\limits_{k=1}^n \frac{ \mathbf{E} \xi_n I_{A_k}}{
\mathbf{E}  \xi_n^2 I_{A_k}}.
\end{eqnarray*}
Theorem \ref{th10} follows from the latter.
$\Box$

Theorem \ref{th10} generalizes the first part of
the classic Borel--Cantelli lemma. If $\{A_n\}$ are independent and
series $\sum\limits_{n=1}^\infty \mathbf{P}(A_n)$ diverges, then
by Theorem \ref{th10}, the upper bound is 1. So, the bound is sharp
in this case.


\smallskip
\noindent
{\bf References}
{\footnotesize

\parindent 0 mm

Andel J., Dupas V., 1989. An extension of the Borel lemma.
Comment.~Math.~Univ.~Carolin. 30, 403–404.


Boros E., Pr\'{e}kopa A., 1989.  Closed form two-sided bounds for probabilities
that at least $r$ and exactly $r$ out of $n$ events occurs.
Math. Oper. Research. 14, 317--342.

de Caen D., 1997. A lower bound on the probability of a union.
Discrete Math. 169, 217--220.

Chandra T.K., 2012. The Borel-Cantelli lemma. Springer, Heidelberg.

Chung K.L., Erd\H{o}s P., 1952. On the application of the Borel-Cantelli lemma.
Trans. Amer. Math. Soc. 72, 179–186.

Dawson D.A., Sankoff D., 1967. An inequality for probabilities.
Proc.~Amer.~Math.~Soc. 18, 504--507.

Erd\H{o}s P., R\'{e}nyi A., 1959. On Cantor’s series with convergent $\sum 1/q$,
Ann. Univ. Sci. Budapest Sect. Math. 2, 93–109.

Feng C., Li  L., Shen  J., 2009. On the Borel--Cantelli lemma and its generalization.
Comptes Rendus Math. 347, 1313–1316.

Frolov A.N., 2012. Bounds for probabilities of unions of events and
the Borel--Cantelli lemma.
Statist. Probab. Lett. 82, 2189--2197.

Frolov A.N., 2014. On inequalities for probabilities of unions of
events and the Borel--Cantelli lemma.
Vestnik Sankt-Peterburgskogo Universiteta, Seriya 1. Matematika, Mekhanika,
Astronomiya, N 2, 21--30. (In Russian) English translation:
Vestnik St.Petersburg University, Mathematics, 2014, N 2, 60--67.
Allerton Press, Inc.

Galambos J., Simonelli I., 1996. Bonferroni-type inequalities with applications.
Springer-Verlag N.Y.

Gallot S. 1966. A bound for the maximum of a number of random variables.
J.~Appl.~Probab. 3, 556--558.

Kochen S., Stone C., 1964. A note on the Borel-Cantelli lemma.
Illinois J. Math. 8, 248–251.

Kounias E.G., 1968. Bounds for the probability of a union, with applications.
Ann. Math. Statist. 39, 2154--2158.

Kwerel S.M., 1975. Bounds on the probability of the union and intersection of
$m$ events.
Adv. Appl. Probab. 7, 431--448.

Kuai H., Alajaji F., Takahara G., 2000. A lower bound on the probability
of a finite union of events.
Discrete Math. 215, 147--158.

M\'{o}ri T.F., Sz\'{e}kely  G.J., 1983.
On the Erd\H{o}s--R\'{e}nyi generalization of the
Borel--Cantelli lemma. Studia Sci. Math. Hungar. 18, 173–182.

Martikainen A.I., Petrov V.V., 1990. On the Borel-Cantelli lemma.
Zapiski Nauch. Semin. LOMI 184, 200--207 (in Russian).
English translation in: J. Math. Sci. 1994, 63, 540--544.


Petrov V.V., 2002. A note on the Borel--Cantelli lemma,
Statist.~Probab.~Lett. 58, 283–286.


Pr\'{e}kopa A., 2009. Inequalities for discrete higher order
convex functions. J.~Math.~Inequalities. 4, 485--498.

Spitzer F., 1964. Principles of random walk. Van Nostrand, Princeton.

Xie Y.Q., 2008. A bilateral inequality for the Borel--Cantelli lemma.
Statist.~Probab.~Lett. 78, 390--395.
}

\end{document}